\documentclass[a4,twoside,10pt]{amsart}

\usepackage{a4wide}
\usepackage{amsfonts}
\usepackage{amsmath}
\usepackage{amssymb}
\usepackage{amsthm}
\usepackage[english]{babel}
\usepackage{enumerate}
\usepackage{fontenc}
\usepackage[colorlinks=true,linkcolor=blue,citecolor=red,urlcolor=blue]{hyperref}
\usepackage[latin1]{inputenc}
\usepackage{orcidlink}
\usepackage{stackrel}
\usepackage[english]{varioref}
\usepackage[all]{xy}
\usepackage{url}

\title{An overview on curve semistable\\ and numerically flat Higgs bundles}

\author{Armando Capasso}
\address{Universit\`a degli Studi di Trieste, P.le Europa 1, Trieste (Italy), C.A.P. 34127}
\email{armando.capasso@units.it}

\thanks{A.C. is member of INdAM-GNSAGA. ORCID: 0009-0001-5463-7221 \orcidlink{0009-0001-5463-7221}}

\date{}

\theoremstyle{plain}
\newtheorem{theorem}{Theorem}[section]

\newtheorem{lemma}[theorem]{Lemma}

\theoremstyle{definition}
\newtheorem{conjecture}{Conjecture}
\newtheorem{definition}[theorem]{Definition}

\newtheorem{ex}[theorem]{Example}
\newenvironment{example}{\begin{ex}}{\hfill{$\triangle$}\end{ex}}
\newtheorem{exs}[theorem]{Examples}
\newenvironment{examples}{\begin{exs}}{\end{exs}}

\newtheorem{rem}[theorem]{Remark}
\newenvironment{remark}{\begin{rem}}{\hfill{$\Diamond$}\end{rem}}
\newtheorem{rems}[theorem]{Remarks}
\newenvironment{remarks}{\begin{rems}}{\end{rems}}
\newtheorem*{Proof}{Proof}
\newenvironment{prf}{\begin{Proof}}{\hfill\text{Q.e.d.}\end{Proof}}

\DeclareMathOperator{\End}{End}
\DeclareMathOperator{\gr}{Gr}
\DeclareMathOperator{\hgr}{\mathfrak{Gr}}
\let\hom\relax
\DeclareMathOperator{\hom}{Hom}
\DeclareMathOperator{\Id}{Id}

\DeclareMathOperator{\rank}{rank}
\DeclareMathOperator{\second}{\prime\prime}

\DeclareMathOperator{\spec}{Spec}

\newcommand{\C}{\mathbb{C}}
\newcommand{\K}{\mathbb{K}}
\newcommand{\N}{\mathbb{N}}
\newcommand{\Q}{\mathbb{Q}}

\newcommand{\Z}{\mathbb{Z}}

\newcommand{\cE}{\mathcal{E}}
\newcommand{\cF}{\mathcal{F}}
\newcommand{\cO}{\mathcal{O}}

\newcommand{\fE}{\mathfrak{E}}
\newcommand{\fF}{\mathfrak{F}}
\newcommand{\fQ}{\mathfrak{Q}}
\newcommand{\fS}{\mathfrak{S}}

\newcommand{\rA}{\mathrm{A}}
\newcommand{\rH}{\mathrm{H}}
\newcommand{\rN}{\mathrm{N}}

\newcommand{\of}{\overline{f}}

\begin{document}

\begin{abstract}
After recalling the basic notions concerning Higgs-Grassmannian schemes, I review how these later can be used to define generalisations of the notion of positivity conditions, such as numerically flatness, which ``feel'' the Higgs field. Then I prove several properties of Higgs bundles, over smooth projective varieties defined over an algebraically closed field of characteristic $0$, satisfying these conditions. Finally, I discuss how one can relate them to semistability of the so-called ``curve semistable'' Higgs bundles.
\end{abstract}

\maketitle

\section*{Introduction}
\markboth{An overview on curve semistable and numerically flat Higgs bundles}{Armando Capasso}

\noindent In order to construct the \emph{moduli space of vector bundles over an irreducible smooth projective curve}, Mumford has introduced the so-called (\emph{slope semi})\emph{stability condition}. Later, Takemoto has extended this notion on any irreducible smooth projective variety (cfr. Definition \ref{def0.1}). In \cite{M:Y}, Miyaoka has related this condition to the \emph{positivity} of a numerical class.\medskip

\noindent To be clear, let $X$ be an irreducible projective variety over an algebraically closed field. A line bundle $L$ over $X$ is \emph{numerically eventually free} (\emph{nef}, for short, see also \cite[Remark 1.4.2]{L:RK}) if for any irreducible curve $C$ on $X$ the following inequality $\displaystyle\int_Cc_1(L)\geq0$ holds. A vector bundle $E$ over $X$ is \emph{nef} if the line bundle $\cO_{\gr_1(E)}(1)$ over $\gr_1(E)$ (the projective bundle of rank $1$ locally free quotients of $E$, see \cite{L:RK}) is nef; and $E$ is \emph{numerically flat} (\emph{nflat}, for short) if $E$ and $E^{\vee}$ are both nef. Thus $E$ is semistable if and only if the \emph{normalized hyperplane class $\lambda_1(E)\in\rN^1\left(\gr_1(E)\right)$} is nef (cfr. Equation \eqref{eq3.1} and \cite[Theorem 3.1]{M:Y}); here $\rN^1\left(\gr_1(E)\right)$ is the \emph{real vector space of real} $1$-\emph{cocycles on} $\gr_1(E)$ \emph{modulo numerical equivalence}.\medskip

\noindent Bruzzo and Hern\'andez Ruip\'erez have generalised in \cite{B:HR} \emph{Miyaoka's criterion} introducing other numerical classes ${\lambda_s(E)\in\rN^1\left(\gr_1\left(Q_{s,E}\right)\right)}$ and $\theta_s(E)\in\rN^1\left(\gr_1(E)\right)$ where $r$ is the rank of $E$, $s\in\{1,\dotsc,r-1\}$ and $Q_{s,E}$ is the \emph{universal rank} $s$ \emph{quotient bundle over} $\gr_s(E)$ (see Equations \eqref{eq3.1} and \eqref{eq3.2}, respectively). They have proved that the semistability of $E$ implies the nefness of all these numerical classes; and, conversely, if at least one of these numerical classes is nef then $E$ is semistable (\cite[Theorem 1.1]{B:HR}). Moreover, they have given a generalisation of this criterion to semistable vector bundles over smooth complex projective varieties.\medskip

\noindent Here I extend this criterion to irreducible smooth projective varieties defined over an algebraically closed field of characteristic $0$ (Theorem \ref{th4.1}). I highlight that the proofs given here are purely algebraic, \emph{i.e.} no analytic method has been used. As a consequence, a vector bundle $E$ over $X$ whose classes $\lambda_s(E)$'s or $\theta_s(E)$'s are nef, respectively, satisfies the following property: the pull-backs of $E$ to any irreducible smooth projective curve are semistable (\emph{curve semistable}, for simplicity; cfr. also Definition \ref{def3.1}). Moreover, a such vector bundle is semistable with respect to all polarizations of $X$ and $c_2(\End(E))=0\in\rA^2(X)$ (the \emph{Abelian group of} $2$-\emph{cocycles on} $X$ \emph{modulo rational equivalence}). Furthermore, this suggests a generalisation to the \emph{Higgs bundles} setting.\medskip

\noindent In simple words: positivity conditions and semistability for vector bundles over irreducible smooth projective varieties are not ``disjoint'' properties. More generally, an analogous phenomena happens in the Higgs bundles framework. In order to generalise all these results, I need to recall main definitions about \emph{Higgs sheaves on projective schemes}.\medskip

\noindent Let $X$ be a projective scheme over an algebraically field of characteristic $0$, let $\Omega^1_X$ be the cotangent sheaf of $X$ and let $\displaystyle\Omega^p_X=\bigwedge^p\Omega^1_X$ for $p\in\N_{\geq1}$.
\begin{definition}
A \emph{Higgs sheaf} $\fE$ on $X$ is a pair $(\cE,\varphi)$ where $\cE$ is an $\cO_X$-coherent sheaf equipped with a morphism $\varphi\colon\cE\to\cE\otimes\Omega^1_X$ called \emph{Higgs field} such that the composition
\begin{displaymath}
\varphi\wedge\varphi\colon\cE\xrightarrow{\varphi}\cE\otimes\Omega^1_X\xrightarrow{\varphi\otimes\Id}\cE\otimes\Omega^1_X\otimes\Omega^1_X\to\cE\otimes\Omega^2_X
\end{displaymath}
vanishes. A Higgs subsheaf of $\fE$ is a $\varphi$-\emph{invariant subsheaf} $\cF$ of $\cE$, \emph{i.e.} $\varphi(\cF)\subseteq\cF\otimes\Omega_X^1$. A \emph{quotient Higgs} of $\fE$ is a quotient sheaf of $\cE$ such that the corresponding kernel is $\varphi$-invariant. A \emph{Higgs bundle} is a Higgs sheaf on $X$ whose underlying coherent sheaf is locally free.
\end{definition}
\noindent Assuming $X$ to be also an irreducible smooth variety, given a Higgs bundle $\fE=(E,\varphi)$ over $X$, Bruzzo and Hern\'andez Ruip\'erez in \cite{B:HR} have introduced closed subschemes $\hgr_t(\fE)$ of $\gr_t(E)$, called the $t$\emph{-th Higgs-Grassmann schemes of} $\fE$, which parametrizes the rank $t$ quotient Higgs bundles of $\fE$. These schemes enjoy a universal property similar to that of the Grassmann bundles. Let $\fQ_{t,\fE}$ be the restriction of the universal rank $t$ quotient bundle of $E$ to $\hgr_t(\fE)$; this is a rank $t$ quotient Higgs bundles of the pullback of $\fE$ over $\hgr_t(\fE)$. They have defined the numerical classes $\lambda_t(\fE)\in\rN^1\left(\gr_1\left(\fQ_{t,\fE}\right)\right)$ and $\theta_t(\fE)\in\rN^1\left(\hgr_1(\fE)\right)$ where $t\in\{1,\dotsc,r-1\}$ (see Equations \eqref{eq3.1} and \eqref{eq3.2}, respectively). So one has an equivalence between the nefness of these numerical classes and the curve semistability of $\fE$ (see Theorems \ref{th3.1} and \ref{th3.2}); even in this setting, $\fE$ is semistable with respect to all polarizations of $X$.\medskip

\noindent To be precise, let $H$ be a polarization of $X$ and let $\fE=(\cE,\varphi)$ be a torsion-free Higgs sheaf on $X$, if not otherwise indicated. One defines the \emph{slope} of $\fE$ as $\displaystyle\mu(\fE)=\frac{\deg(\fE)}{\rank(\fE)}\in\Q$, where $\displaystyle\deg(\fE)=\int_Xc_1(\cE)\cdot H^{n-1}\in\Z$. In a similar way as for vector bundles there is a notion of stability for torsion-free Higgs sheaves.
\begin{definition}\label{def0.1}
$\fE$ is \emph{semistable} (respectively, \emph{stable}) if $\mu(\fE)\stackrel[(<)]{\textstyle\leq}{}\mu(\fQ)$ for every torsion-free quotient Higgs sheaf $\fQ$ with $0<\rank(\fQ)<\rank(\fE)$. In the other eventuality, $\fE$ is \emph{unstable}.
\end{definition}
\noindent On another hand, these Higgs-Grassmann schemes of a Higgs bundle $\fE=(E,\varphi)$ are an ingredient of the positivity conditions for Higgs bundles. The idea is the following: $\det(E)$ has to be nef. Of course this is not enough if $r\geq2$, so one requires also the ``\emph{Higgs nefness}'' of all \emph{universal quotient Higgs bundles} $\fQ_{t,\fE}$ recursively (see Definition \ref{def1.1}). A such Higgs bundle $\fE$ is called \emph{Higgs nef} (\emph{H-nef}, for short); and it follows to define \emph{Higgs numerically flat} (\emph{H-nflat}, for short) a H-nef Higgs bundle $\fE$ whose dual Higgs bundle $\fE^{\vee}$ is H-nef as well. This idea has been developed by Bruzzo, Hern\'andez Ruip\'erez and Gra\~na Otero in \cite{B:HR,B:GO:1}.\medskip

\noindent The interplay between semistability and positivity conditions for Higgs bundles has been studied by Biswas, Bruzzo, myself, Gra\~na Otero, Gurjar, Hern\'andez Ruip\'erez, Lanza, Lo Giudice and Peragine in \cite{B:HR,B:GO:1,B:GO:2,B:GO:3,B:LG,L:LG,B:B:G,B:L:LG,B:P,B:C,B:GO:HR,B:C:GO}. In particular, it is open the problem whether a curve semistable Higgs bundle $\fE=(E,\varphi)$ has $c_2(\End(E))=0\in\rA^2(X)$ (Theorem \ref{th4.3}). Equivalently, it is an open problem the vanishing of Chern classes of H-nflat Higgs bundles (\cite[Corollary 3.2]{B:B:G}). For a recap on this open problem, I remind to Remark \ref{rem4.1}.\medskip

\noindent An equivalent way to define H-nflat Higgs bundles is the following one: a Higgs bundle $\fE=(E,\varphi)$ over an irreducible smooth projective variety is H-nflat if it is curve semistable and $c_1(E)\equiv_{num}0$ (cfr. Lemma \ref{lem2.2}.\ref{lem2.2.a}). So this combined with \cite[Corollary]{L:A} simplifies the problem to prove the vanishing of $c_2(E)\in\rA^2(X)$. Other simplifications of this problem are recalled and proved at the end of this paper.\medskip

\noindent \textbf{Notations and conventions.} $\K$ is an algebraically closed field of characteristic $0$, unless otherwise indicated. By a projective variety $X$ I mean a projective integral scheme over $\K$ of dimension $n\geq1$ and of finite type. If $n\in\{1,2\}$ I shall write projective curve or projective surface, respectively. Whenever I consider a morphism $f\colon C\to X$, I understand that $C$ is an irreducible smooth projective curve. $\rN^1(X)$ is the real vector space of real $1$-cocycles on $X$ modulo numerical equivalence, $\rA^2(X)$ is the Abelian group of $2$-cocycles on $X$ modulo rational equivalence.\medskip

\noindent \textbf{Acknowledgements.} The main results of this paper are based on my Ph.D. thesis \cite{AC:PhD} written under supervision of U. Bruzzo and B. Gra\~na Otero. I am very grateful to them for their help, their energy and their support. The text of this paper reflects the mini talk I gave on August 8th 2025 during the ``Workshop on Character Varieties and Higgs Bundles'' held at Universidad de Costa Rica. I am grateful to the organizers for their invitation.\medskip

\noindent \textbf{Further informations.} A.C. is member of INdAM-GNSAGA. ORCID: 0009-0001-5463-7221 \orcidlink{0009-0001-5463-7221}\medskip

\noindent{\bf Statement about competing or financial interests.} The author has no competing or financial interests to declare that are relevant to the content of this article.

\section{Higgs-Grassmann schemes and\\ H-ample, H-nef and H-nflat Higgs bundles}
\markboth{An overview on curve semistable and numerically flat Higgs bundles}{Armando Capasso}

\noindent Let $\fE=(E,\varphi)$ be a rank $r\geq2$ Higgs bundle over a smooth projective variety $X$, and let ${s\in\{1,\dotsc,r-1\}}$ be an integer number. Let $p_s\colon\gr_s(E)\to X$ be the \emph{Grassmann bundle} parametrizing rank $s$ locally free quotients of $E$ (see \cite{L:RK}). Consider the short exact sequence of vector bundles over $\gr_s(E)$ 
\begin{displaymath}
\xymatrix{
0\ar[r] & S_{r-s,E}\ar[r]^(.55){\eta} & p_s^{\ast}E\ar[r]^(.45){\epsilon} & Q_{s,E}\ar[r] & 0
},
\end{displaymath}
where $S_{r-s,E}$ is the \emph{universal rank $r-s$ subbundle} and $Q_{s,E}$ is the \emph{universal rank $s$ quotient bundle of} $p_s^{\ast}E$, respectively. One defines the closed subschemes $\hgr_s(\fE)\subseteq\gr_s(E)$ (the $s$\emph{-th Higgs-Grassmann schemes of} $\fE$) as the zero loci of the composite morphisms
\begin{displaymath}
(\epsilon\otimes\Id)\circ p_s^{\ast}\varphi\circ\eta\colon S_{r-s,E}\to Q_{s,E}\otimes p_s^{\ast}\Omega_X^1.
\end{displaymath}
The restriction of the previous sequence to $\hgr_s(\fE)$ yields a universal short exact sequence
\begin{displaymath}
\xymatrix{
0\ar[r] & \mathfrak{S}_{r-s,\fE}\ar[r]^(.55){\psi} & \rho_s^{\ast}\fE\ar[r]^(.45){\eta} & \fQ_{s,\fE}\ar[r] & 0,
}
\end{displaymath}
where $\fQ_{s,\fE}=Q_{s,E|\hgr_s(\fE)}$ is equipped with the quotient Higgs field induced by $\rho_s^{\ast}\varphi$, and ${\rho_s=p_{s|\hgr_s(\fE)}}$. The scheme $\hgr_s(\fE)$ enjoys the usual universal property: for a morphism of varieties $f\colon Y\to X$, the morphism $g\colon Y \to \gr_s(E)$ given by a rank $s$ quotient $Q$ of $f^{\ast}E$ factors through $\hgr_s(\fE)$ if and only if $\varphi$ induces a Higgs field on $Q$.
\begin{definition}[{see \cite[Definition A.2]{B:GO:3}}]\label{def1.1}
A Higgs bundle $\fE=(E,\varphi)$ of rank one is said to be \emph{Higgs-ample}/\emph{Higgs-numerically effective} (\emph{H-ample}/\emph{H-nef}, for short) if $E$ is \emph{ample}/\emph{numerically effective} in the usual sense. If $\rank(\fE)\geq 2$, we inductively define H-ampleness/H-nefness by requiring that
\begin{enumerate}[a)]
\item the determinant line bundle $\det(E)$ is ample/nef, and
\item all Higgs bundles $\fQ_{s,\fE}$ are H-ample/H-nef for all $s$.
\end{enumerate}
$\fE$ is \emph{Higgs-numerically} flat (\emph{H-nflat}, for short) if $\fE$ and $\fE^{\vee}$ are both H-nef.
\end{definition}
\begin{remarks}\label{rem1.1}
\,\begin{enumerate}[a)]
\item Note that if $E$ nef/nflat in the usual sense, then $\fE$ is H-nef/H-nflat. If $\varphi=0$, the Higgs bundle $\fE=(E,0)$ is H-nef/H-nflat if $E$ is nef/nflat in the usual sense.
\item\label{rem1.1.b} The recursive condition in this definition can be recasted as follows. Let ${1\leq s_1<s_2<\dotsc<s_k<r}$ and let $\fQ_{\left(s_1 \cdots,s_k\right),\fE}$ be the rank $s_1$ universal quotient Higgs bundle over $\hgr_{s_1}\left(\fQ_{\left(s_2,\dotsc,s_k\right),\fE}\right)$, obtained by taking the successive universal quotient Higgs bundles of $\fE$ of rank $s_k$, then $s_{k-1}$, all the way to rank $s_1$. The H-ampleness/H-nefness condition for $\fE$ is equivalent to requiring that all line bundles $\det(\fE)$ and $\det(\fQ_{\left(s_1,\cdots,s_k\right),\fE})$ are ample/nef.
\item\label{rem1.1.c} The first Chern class of an H-nflat Higgs bundle is numerically zero, because the corresponding determinant bundle is nflat. \hfill{$\Diamond$}
\end{enumerate}
\end{remarks}
\noindent I refer the reader to \cite[Remark 4.2]{B:GO:2}, \cite[Example 3.3]{B:C:GO}, examples \ref{ex2.1} and \ref{ex3.1} for furthermore considerations on the previous definitions.

\section{Properties of H-nef and H-nflat Higgs bundles}
\markboth{An overview on curve semistable and numerically flat Higgs bundles}{Armando Capasso}

\noindent Higgs bundles over $X$ which are H-nef satisfy properties analogous to those of nef vector bundles. These properties have been proved in \cite{B:GO:1,B:GO:2,B:B:G,B:C:GO}; here I list some of them for completeness.
\begin{lemma}\label{lem2.1}
Let $\fE=(E,\varphi)$ be an H-nef Higgs bundle over $X$. The following statements hold.
\begin{enumerate}[a)]
\item\label{lem2.1.a} Let $f\colon Y\to X$ be a morphism of smooth projective varieties. Then $f^{\ast}\fE$ is H-nef (\cite[Proposition 2.6.(ii)]{B:GO:1}). If $f$ is also surjective and $f^{\ast}\fE$ is H-nef then $\fE$ is H-nef (\cite[Lemma 3.4]{B:B:G}).
\item\label{lem2.1.b} Every quotient Higgs bundle of $\fE$ is H-nef (\cite[Lemma 3.5]{B:B:G}).
\item Tensor products of H-nef Higgs bundles are H-nef (\cite[Theorem 3.6]{B:B:G}).
\item Exterior and symmetric powers of H-nef Higgs bundles are H-nef (cfr. \cite[Propositions 3.5, 4.4 and Lemma 4.5]{B:GO:2}).
\item\label{lem2.1.e} $\fE$ is H-nef if and only if the Higgs bundle $\fE\otimes\cO_X(D)=(E\otimes\cO_X(D),\varphi\otimes Id)$ is H-ample for every ample Cartier $\Q$-divisor $D$ in $X$ (\cite[Proposition 2.6.(i)]{B:GO:1}).
\item $\fE$ is H-nef if and only if $\det(E)$ is nef and for every finite morphism $f\colon C\to X$, any every Higgs quotient bundle $f^{\ast}\fE\to\fQ\to0$, the inequality $\deg\fQ\geq0$ holds (\cite[Corollary 3.7]{B:C:GO}).
\item The extensions of H-nef Higgs bundles are H-nef (cfr. \cite[Theorem 3.9]{B:C:GO}).
\end{enumerate}
\end{lemma}
\begin{remark}
In \cite{B:GO:1}, the first part of Lemmata \ref{lem2.1}.\ref{lem2.1.a} and \ref{lem2.1}.\ref{lem2.1.e} are proved assuming that $\K=\C$, however these proofs work in general by \cite[Propositions 6.1.2.(ii) and 6.1.8.(iv)]{L:RK}.
\end{remark}
\noindent Category of H-nflat Higgs bundles satisfies other properties which have been proved in \cite{B:GO:3,B:B:G,L:LG,B:C}.
\begin{lemma}\label{lem2.2}
\,
\begin{enumerate}[a)]
\item\label{lem2.2.a} $\fE$ is H-nflat if and only if the pullback of $\fE$ via any $f\colon C\to X$ is semistable and $\displaystyle\int_Cf^{\ast}c_1(E)=0$ (cfr. \cite[Lemma A.7]{B:GO:3}).
\item\label{lem2.2.b} Any H-nflat Higgs bundle is semistable (\cite[Proposition A.8]{B:GO:3}).
\item\label{lem2.2.c} Extensions of H-nflat Higgs bundles are H-nflat (\cite[Proposition 3.1.(iii)]{B:B:G}).
\item\label{lem2.2.d} Tensor products of H-nflat Higgs bundles are H-nflat (\cite[Proposition 3.1.(iv)]{B:B:G}).
\item Kernels and cokernels of morphisms of H-nflat Higgs bundles are H-nflat Higgs bundles (\cite[Propositions 3.7 and 3.8]{B:B:G}).
\item\label{lem2.2.f} Let $\fE$ be a Higgs bundle over $X$. $\fE$ is H-nflat if and only if it is \emph{pseudostable}, \emph{i.e.} it has a filtration
\begin{equation}\label{eq2.1}
0=\fF_{m+1}\subsetneqq\fF_m\subsetneqq\dotsc\subsetneqq\fF_0=\fE
\end{equation}
whose quotients $\fF_i/\fF_{i+1}$ are locally free and stable, and these quotients are also H-nflat (\cite[Theorem 3.2]{B:C}).
\end{enumerate}
\end{lemma}
\begin{remark}
The ``if part'' of Lemma \ref{lem2.2}.\ref{lem2.2.f} follows by Lemma \ref{lem2.2}.\ref{lem2.2.c}. Indeed, let $\fE$ be a Higgs bundle and let
\begin{displaymath}
0=\fF_0\subsetneqq\fF\subsetneqq\fF_1\subsetneqq\dotsc\subsetneqq\fF_m\subsetneqq\fF_{m+1}=\fE
\end{displaymath}
be a filtration of $\fE$ in Higgs subbundle whose quotients $\fF,\fQ_1,\dotsc,\fQ_{m+1}$ are locally free, stable and H-nflat. As explained in the proof of \cite[Theorem 3.2]{B:C}, $\fF_1$ is H-nflat. Consider the short exact sequence
\begin{displaymath}
\xymatrix{
0\ar[r] & \fF_1\ar[r] & \fF_2\ar[r] & \fQ_2\ar[r] & 0
},
\end{displaymath}
since $\fF_1$ and $\fQ_2$ are H-nflat, by Lemma \ref{lem2.2}.\ref{lem2.2.c} $\fF_2$ is H-nflat. Iterating this reasoning, one proves that $\fE$ is H-nflat.
\end{remark}
\noindent Even if \cite[Lemma A.7 and Proposition A.8]{B:GO:3} and \cite[Theorem 3.2]{B:C} have been proved where $\K=\C$, I shall explain that this hypothesis is unnecessary, in the sense these statements work on algebraically closed fields of characteristic $0$. Thus, also the Lemmata \ref{lem2.2}.\ref{lem2.2.c} and \ref{lem2.2}.\ref{lem2.2.d} hold on any algebraically closed field of characteristic $0$. The original proofs of the remaining lemmata continue to be valid in this extending setting, therefore I do not repeat them here.
\begin{examples}\label{ex2.1}
Let $C$ be a smooth projective curve of genus $g\geq2$.
\begin{enumerate}[a)]
\item\,(cfr. \cite[Example 3.2.9]{AC:PhD}) Let $\fE=(E=L_1\oplus L_2,\varphi)$ be a rank $2$ Higgs bundle over $C$ such that $\varphi(L_1)\subseteq L_2\otimes\Omega^1_X$ and $\varphi(L_2)=0$. By \cite[Proposition 3.2.2]{AC:PhD}, $\fE$ has only a quotient Higgs bundle which is $(L_1,0)$. If one takes $g=2$, $\deg\left(L_1\right)=1$ and $\deg\left(L_2\right)=0$, then $\fE$ is a H-ample Higgs bundle whose underlying vector bundle is not ample.
\item\,(see \cite[Example 2.9]{B:LG}) Let $\fE=\left(E=K^{\frac{1}{2}}\oplus K^{-\frac{1}{2}},\varphi\right)$, where $K$ is the canonical bundle of $C$, $K^{\frac{1}{2}}$ is a line bundle over $C$ whose square is $K$ and
\begin{displaymath}
\varphi=\begin{pmatrix}
0 & \omega\\
1 & 0
\end{pmatrix},\,1\in\hom\left(K^{\frac{1}{2}},K^{-\frac{1}{2}}\otimes K\right),\,\omega\in\rH^0(C,K^2).
\end{displaymath}
$\fE$ is a stable Higgs bundle, because there are no subbundles of positive degree preserved by $\varphi$. Indeed, let $L$ be a line subbundle of $E$ then $L$ is either $K^{-\frac{1}{2}}$ or $K^{\frac{1}{2}}$. Since
\begin{displaymath}
\deg\left(K^{\frac{1}{2}}\right)=g-1>0,\,\deg\left(K^{-\frac{1}{2}}\right)=1-g<0
\end{displaymath}
$L$ destabilizes $E$ if and only if $L=K^{\frac{1}{2}}$. However $K^{\frac{1}{2}}$ does not destabilize $\fE$, because this is not $\varphi$-invariant. Furthermore, since $\deg(\fE)=0$ then $\fE$ is H-nflat by Lemma \ref{lem2.2}.\ref{lem2.2.a}; but the underlying vector bundle is unstable hence it is not nflat (see lemma \ref{lem2.2}.\ref{lem2.2.b}). \hfill{$\triangle$}
\end{enumerate}
\end{examples}
\noindent The category of H-ample Higgs bundles over $X$ has been studied in \cite{B:C:GO}, hence I refer the reader to there for any further detail.

\section{Curve semistable Higgs bundles}
\markboth{An overview on curve semistable and numerically flat Higgs bundles}{Armando Capasso}

\noindent From now on, let $\fE=(E,\varphi)$ be a rank $r$ Higgs bundle over a smooth projective polarized variety $(X,H)$, and let $\cE$ be the sheaf of sections of $E$, if not otherwise indicated.\medskip

\noindent Lemma \ref{lem2.2}.\ref{lem2.2.a} justifies the next definition which specializes Definition \ref{def0.1}.
\begin{definition}\label{def3.1}
$\fE$ is \emph{curve semistable} if for every morphism $f\colon C\to X$ the pullback Higgs bundle $f^{\ast}\fE$ is semistable.
\end{definition}
\noindent The study of this class of Higgs bundle has started in \cite{B:HR} as a generalisation to Higgs bundles framework of Miyaoka's work \cite[Section 3]{M:Y} on semistable vector bundles over smooth projective curves. There Bruzzo and Hern\'andez Ruip\'erez have introduced the following numerical classes
\begin{gather}
\lambda_s(\fE)=\left[c_1\left(\cO_{\gr_1\left(Q_{s,\fE}\right)}(1)\right)\right]-\frac{1}{r}\varpi_s^{\ast}(c_1(E))\in\rN^1\left(\gr_1\left(Q_{s,\fE}\right)\right)\label{eq3.1}\\
\theta_s(\fE)=\left[c_1\left(Q_{s,\fE}\right)\right]-\frac{s}{r}\rho_s^{\ast}(c_1(E))\in\rN^1\left(\hgr_s(\fE)\right)\label{eq3.2}
\end{gather}
where $s\in\{1,\dotsc,r-1\}$ and $\varpi_s\colon\gr_1\left(Q_{s,\fE}\right)\to\hgr_s(\fE)\xrightarrow{\rho_s}X$.\medskip

\noindent I recall that a class $\gamma\in\rN^1(X)$ is \emph{numerically effective} (\emph{nef}, for short) if for any irreducible projective curve $C\subseteq X$ the inequality $\gamma\cdot[C]\geq0$ holds. And one calls $\gamma$ \emph{strictly nef} class if for any irreducible projective curve $C\subseteq X$ the inequality $\gamma\cdot[C]>0$ holds. Bruzzo, Gra\~na Otero and Hern\'andez Ruiper\'ez have proved in \cite{B:HR,B:GO:1,B:GO:3} the following theorems in the complex setting, here I extend these results to smooth projective varieties defined over an algebraically closed field of characteristic $0$, without using any analytic method.
\begin{theorem}[{see \cite[Theorem 1.2]{B:HR}}]\label{th3.1}
$\fE$ is curve semistable if and only if each class $\theta_s(\fE)$ is \emph{nef}, where $s\in\{1,\dotsc,r-1\}$.
\end{theorem}
\begin{theorem}[{see \cite[Theorem 1.2]{B:HR}}]\label{th3.2}
$\fE$ is curve semistable if and only if each class $\lambda_s(\fE)$ is \emph{nef}, where $s\in\{1,\dotsc,r-1\}$.
\end{theorem}
\begin{remark}
By definition $\lambda_1(\fE)=\theta_1(\fE)$. Thus if $r=2$ the previous theorems are the same.
\end{remark}
\noindent Before proving these theorems, I recall the following lemma.
\begin{lemma}[{see \cite[Lemma 3.3]{B:HR}}]\label{lem3.1}
Let $f\colon Y\to X$ be a finite surjective morphism of smooth projective curves. Then $\fE$ is semistable if and only if $f^{\ast}\fE$ is semistable.
\end{lemma}
\begin{prf}
If $\fE$ is unstable then there exists a torsion-free Higgs subsheaf $\fF$ of $\fE$ such that ${\mu(\fF)>\mu(\fE)}$, hence $\mu\left(f^{\ast}\fF\right)>\mu\left(f^{\ast}\fE\right)$ \emph{i.e.} $f^{\ast}\fE$ is unstable. Let me assume $f^{\ast}\fE$ unstable, without loss of generality one can assume that $f$ is a Galois covering\footnote{
In general, if $char(\K)>0$ then one needs also that $f$ is separable as stated in \cite{M:Y,B:HR}.
} with Galois group $G$, \emph{i.e.} the field extension $f^{\#}\colon\K(X)\to\K(Y)$ is normal (and separable) and the relevant Galois group is $G$. Let $\fF=\left(F,f^{\ast}\varphi_{\vert F}\right)$ be the maximal destabilizing Higgs subsheaf of $f^{\ast}\fE$. For any $g\in G$, $g^{\ast}\fF$ is a destabilizing Higgs subsheaf of $f^{\ast}\fE$ of maximal rank; however, by the unicity of the HN-filtration of $f^{\ast}\fE$, it has to be $g^{\ast}\fF=\fF$. From all this, it follows that $F=f^{\ast}E_0$ for some destabilizing subbundle $E_0$ of $E$. By \cite[Exercise III.9.3.a]{H:RC} $f$ is a flat morphism, and since it is also surjective then $f$ is (by definition) faithfully flat. Thus the composition $E_0\otimes\Omega^1_X\to E\otimes\Omega^1_X\to(E/E_0)\otimes\Omega^1_X$ vanishes if and only if the composition $F\otimes f^{\ast}\Omega^1_X\to f^{\ast}E\otimes f^{\ast}\Omega^1_X\to(f^{\ast}E/F)\otimes f^{\ast}\Omega^1_X$ vanishes. Consider the following diagram
\begin{displaymath}
\xymatrix{
& & 0\ar[d] & 0\ar[d]\\
& F\ar[r]\ar[dr]_{f^{\ast}\varphi} & f^{\ast}E\otimes f^{\ast}\Omega^1_X\ar[r]\ar[d] & (f^{\ast}E/F)\otimes f^{\ast}\Omega^1_X\ar[d]\\
0\ar[r] & F\otimes\Omega^1_Y\ar[r] & f^{\ast}E\otimes\Omega^1_Y\ar[r] & (f^{\ast}E/F)\otimes\Omega^1_Y\ar[r] & 0
}
\end{displaymath}
since $f^{\ast}\varphi_{\vert F}$ takes values in $F\otimes\Omega^1_Y$ one has the claim, \emph{i.e.} $\fE$ is unstable.
\end{prf}
\noindent \textbf{Proof of Theorem \ref{th3.1}.} Let $C\subseteq X$ be an irreducible projective curve and assume that the restriction $\theta_s(\fE)_{\vert C}$ of $\theta_s(\fE)$ to $\hgr_s\left(\fE_{\vert C}\right)$ is nef. If $\fF=(F,\psi)$ is a rank $s$ torsion-free Higgs quotient sheaf of $\fE_{\vert C}$ then by \cite[Corollary at page 75]{O:S:S} $F$ is locally free. By Universal Property of Higgs-Grassmann schemes, there exists a unique section $\sigma\colon C\to\hgr_s\left(\fE_{\vert C}\right)$ such that $\fF=\sigma^{\ast}\fQ_{s,\fE\vert C}$, where $\fQ_{s,\fE\vert C}$ is the restriction of $\fQ_{s,\fE}$ to $\hgr_s\left(\fE_{\vert C}\right)$. Then
\begin{displaymath}
0\leq\theta_s(\fE)_{\vert C}\cdot[\sigma(C)]=\deg(F)-\frac{s}{r}\deg\left(E_{\vert C}\right)=s\left(\mu(F)-\mu\left(E_{\vert C}\right)\right);
\end{displaymath}
thus if each $\theta_s(\fE)$ is nef then $\fE_{\vert C}$ is semistable, \emph{i.e.} $\fE$ is curve semistable. \emph{Vice versa}, let $\fE$ be curve semistable and let me assume that $\theta_s(\fE)$ is not nef for some $s\in\{1,\dotsc,r-1\}$ \emph{i.e.} there exists an irreducible projective curve $C^{\prime}\subseteq\hgr_s(\fE)$ such that $\theta_s(\fE)\cdot\left[C^{\prime}\right]<0$. Under this hypothesis, $C^{\prime}$ is not contained in a fibre of $\rho_s$, hence it surjects onto a projective curve $C\subseteq X$. One may choose a projective curve $C^{\second}$ and a morphism $h\colon C^{\second}\to C$ such that $\widetilde{C}=C^{\second}\times_CC^{\prime}$ is a union of projective curves $C_j$ isomorphic to $C$ (see the proof of \cite[Theorem 3.1]{M:Y}); $\theta_s\left(h^{\ast}\fE\right)\cdot[C_j]<0$ for any index $j$ evidently. For clarity, one has the following Cartesian diagram
\begin{displaymath}
\xymatrix{
\hgr_s\left(h^{\ast}\fE\right)\ar[r]_{h_s}\ar[d] & \hgr_s\left(\fE_{\vert C}\right)\ar@{^{(}->}[r]\ar[d] & \hgr_s(\fE)\ar[d]^{\rho_s} & C^{\prime}\ar@{_{(}->}[l]\\
C^{\second}\ar[r]_{h} & C\ar@{^{(}->}[r]_{i} & X
},
\end{displaymath}
by the universal property of fibre products $\widetilde{C}\hookrightarrow\hgr_s\left(h^{\ast}\fE\right)$. Let $\fE_j=h_s^{\ast}\left(\rho_s^{\ast}\fE_{\vert C}\right)_{\vert C_j}$ and let $\fQ_j=\fQ_{s,h^{\ast}\fE\vert C_j}$. By Lemma \ref{lem3.1} $\fE_j$ is a semistable Higgs bundle and in particular $\mu(\fQ_j)\geq\mu(\fE_j)$. On the other hand
\begin{displaymath}
0>\theta_s\left(h^{\ast}\fE\right)\cdot[C_j]=\left(c_1(Q_j)-\frac{s}{r}\rho_s^{\ast}c_1(E)\right)\cdot[C_j]=s\left(\mu(Q_j)-\mu(E_j)\right)
\end{displaymath}
where $E_j$ and $Q_j$ are the underlying vector bundles to $\fE_j$ and $\fQ_j$, respectively. But this is a contradiction with the assumptions, hence all classes $\theta_s(\fE)$ have to be nef. \hfill\text{Q.e.d.}\medskip

\noindent \textbf{Proof of Theorem \ref{th3.2}.} The class $\lambda_s(\fE)$ may be regarded as the numerical class of the hyperplane bundle of the Higgs $\Q$-bundle $\fF_s=\fQ_{s,\fE}\otimes\rho_s^{\ast}\left(\det(E)^{-1/r}\right)$ over $\hgr_s(\fE)$. As a consequence
\begin{displaymath}
c_1(\fF_s)=c_1\left(Q_{s,\fE}\right)-\frac{s}{r}\rho_s^{\ast}c_1(E)\Rightarrow\theta_s(\fE)=\left[c_1(\fF_s)\right]\in\rN^1\left(\hgr_s(\fE)\right).
\end{displaymath}
From all this, the classes $\lambda_s(\fE)$ are nef if and only if the classes $\theta_s(\fE)$ are nef, \emph{i.e.} if and only if $\fE$ is curve semistable (Theorem \ref{th3.1}). \hfill\text{Q.e.d.}
\begin{remark}
Mimicking Definition \ref{def3.1}, $\fE$ is \emph{curve stable} if for every morphism $f\colon C\to X$ the pullback Higgs bundle $f^{\ast}\fE$ is stable.\smallskip

\noindent By the previous proofs, the curve stability of $\fE$ implies the strictly nefness of $\theta_s(\fE)$ and $\lambda_s(\fE)$ for any $s\in\{1,\dotsc,r-1\}$.
\end{remark}
\noindent \textbf{Proof of Lemma \ref{lem2.2}.\ref{lem2.2.a}}
Let assume me $\fE$ satisfies the hypotheses, the classes $\displaystyle c_1\left(\cO_{\gr_1\left(Q_{s,\fE}\right)}(1)\right)$ are nef for any $s\in\{1,\dotsc,r-1\}$ (Theorem \ref{th3.2}), hence their pullbacks to $\hgr_s\left(f^{\ast}\fQ_{s,\fE}\right)$ are nef. In other words, using the notations introduced in Remark \ref{rem1.1}.\ref{rem1.1.b}, $\fQ_{(1,s),\fE}$ is nef; and it remains only to prove that $\det\left(\fQ_{{(s_1,\cdots, s_k)},\fE}\right)$ are nef for all strings of integer numbers $1\leq s_1<s_2<\dotsc<s_k<r$. By construction, $\fQ_{(s_1,\cdots, s_k),\fE}$ is a Higgs bundle over $\hgr_{s_1}\left(\fQ_{(s_2,\cdots, s_k),\fE}\right)$ and there is a morphism $\rho_{(s_1,\cdots, s_k)}\colon\hgr_{s_1}\left(\fQ_{(s_2,\cdots, s_k),\fE}\right)\to X$ such that $\fQ_{(s_1,\cdots, s_k),\fE}$ is a rank $s_1$ quotient Higgs bundle of $\rho_{(s_1,\cdots, s_k)}^{\ast}\fE$. Thus there exists a unique morphism $g_{(s_1,\cdots, s_k)}\colon\hgr_{s_1}\left(\fQ_{(s_2,\cdots, s_k),\fE}\right)\to\hgr_{s_1}(\fE)$ such that $\fQ_{(s_1,\cdots, s_k),\fE}=g_{(s_1,\cdots, s_k)}^{\ast}\fQ_{s_1,\fE}$. Therefore
\begin{displaymath}
\left[c_1\left(\fQ_{(s_1,\cdots, s_k),\fE}\right)\right]=g_{(s_1,\cdots, s_k)}^{\ast}\left[c_1\left(\fQ_{s_1,\fE}\right)\right]=g_{(s_1,\cdots, s_k)}^{\ast}\theta_{s_1}(\fE)
\end{displaymath}
because $\displaystyle\int_Cf^{\ast}c_1(E)=0$, hence $\det\left(\fQ_{(s_1,\cdots, s_k),\fE}\right)$ is nef because $\theta_{s_1}(\fE)$ is nef by Theorem \ref{th3.1}. In other words $\fE$ is H-nef. Repeating all this reasoning, considering that also $\fE^{\vee}$ is curve semistable, $c_1\left(E^{\vee}\right)=-c_1(E)$ hence $\displaystyle\int_Cf^{\ast}c_1\left(E^{\vee}\right)=0$, one proves in the same way the H-nefness of $\fE^{\vee}$, \emph{i.e.} $\fE$ is H-nflat.\smallskip

\noindent \emph{Vice versa} holds, since the proof of \cite[Proposition A.8]{B:GO:3} works on smooth projective varieties defined over an algebraically closed field of characteristic $0$. \hfill Q.e.d.
\begin{example}\label{ex3.1}
Let $X$ be a minimal smooth surface such that the canonical bundle $K_X$ is nef as well as being big, \emph{i.e.} $X$ is also a surface of general type. The Chern classes of $X$ satisfy the \emph{Bogomolov-Miyaoka-Yau inequality} (\cite[Proposition 7.1]{M:Y})
\begin{equation*}
BMY(X)\stackrel{def.}{=}\int_X3c_2(X)-c_1(X)^2\geq0,
\end{equation*}
where $c_k(X)\stackrel{def.}{=}(-1)^kc_k\left(\Omega^1_X\right)$ for any $k$. One considers the so-called \emph{Simpson system}, \emph{i.e.} the Higgs bundle $\fS=(S,\varphi)$, where $S=\Omega^1_X\oplus\cO_X$ and
\begin{displaymath}
\varphi=\begin{pmatrix}
0 & 0\\
\Id & 0
\end{pmatrix},\,\Id\in\hom\left(\Omega^1_X,\Omega^1_X\right).
\end{displaymath}
If $BMY(X)=0$ then the Higgs bundle $\fS(-\beta K_X)$ is H-nef for every rational number\footnote{
As it is customary, I formally consider twistings by rational divisors, which make sense after pulling back to a (possibly ramified) finite covering of $X$; on the other hand, the properties of being semistable, H-nef are invariant under such coverings (see Lemmata \ref{lem2.1}.\ref{lem2.1.a} and \ref{lem3.1}).
} $\displaystyle0<\beta\leq\frac{1}{3}$, but $S\left(-\beta K_X\right)$ is not nef (\cite[Theorem 4.2]{B:C:GO}). Moreover, $\End(\fS)$ is H-nflat, by \cite[Proposition 4.1]{B:C:GO} and Lemma \ref{lem2.2}.\ref{lem2.2.a}, but $\End(S)$ is not nflat since $S$ is not semistable (see Lemma \ref{lem2.2}.\ref{lem2.2.b}).
\end{example}
\begin{remark}
In \cite{A:C:2}, last example has been extended to $n$-dimensional varieties $X$ whose canonical bundle is ample and $2(n+1)c_2(X)=nc_1(X)^2\in\rA^2(X)$ (cfr. \cite[Corollary 3.4.b]{A:C:2}).
\end{remark}

\section{Curve semistable Higgs bundles whose discriminant class vanishes}
\markboth{An overview on curve semistable and numerically flat Higgs bundles}{Armando Capasso}

\noindent I have recalled the definition of curve semistable Higgs bundles, because I would like to extend the following theorem to Higgs bundles setting.
\begin{theorem}[{cfr. \cite[Theorem 2]{N:Nak} and \cite[Theorem 1.4]{B:HR}}]\label{th4.1}
For a vector bundle $E$ over a smooth projective variety $X$ the following statements are equivalent:
\begin{enumerate}[a)]
\item $\theta_1(E)$ is nef;
\item\label{th4.1.b} $E$ is curve semistable;
\item\label{th4.1.c} $E$ is semistable with respect to some polarization $H$ and $c_2\left(\End(E)\right)=0\in\rA^2(X)$;
\item\label{th4.1.d} $E$ is semistable with respect to some polarization $H$ and $\displaystyle\int_Xc_2\left(\End(E)\right)\cdot H^{n-2}=0$.
\end{enumerate}
\end{theorem}
\begin{prf}
\textbf{(a) is equivalent to (b).} This is \cite[Theorem 3.1]{M:Y} when $\dim X=1$. Let ${\dim X\geq2}$. If $E$ is curve semistable then $\theta_1(E)$ is nef by Theorem \ref{th3.1}. \emph{Vice versa}, if $\theta_1(E)$ is nef, let $f\colon C\to X$ be a morphism. Consider the following Cartesian diagram
\begin{displaymath}
\xymatrix{
\gr_1\left(f^{\ast}E\right)\ar[r]^{\of}\ar[d]_{p_{1\vert C}\equiv f^{\ast}\left(p_1\right)} & \gr_1(E)\ar[d]^{p_1}\\
C\ar[r]_f & X
},
\end{displaymath}
one has
\begin{displaymath}
\of^{\ast}\lambda_1(E)=\of^{\ast}\left(\left[c_1\left(Q_{1,E}\otimes p_1^{\ast}\left(\det(E)^{-1/r}\right)\right)\right]\right)=\left[c_1\left(Q_{1,f^{\ast}E}\otimes p_{1\vert C}^{\ast}\left(\det(E)^{-1/r}\right)\right)\right]=\lambda_1\left(f^{\ast}E\right).
\end{displaymath}
Since $Q_{1,E}\otimes p_1^{\ast}\left(\det(E)^{-1/r}\right)$ is nef then its pullback $Q_{1,f^{\ast}E}\otimes p_{1\vert C}^{\ast}\left(\det(E)^{-1/r}\right)$ via $f$ is nef (\cite[Proposition 2.6.(ii)]{B:GO:1}). Thus $\lambda_1\left(f^{\ast}E\right)$ is nef as well, and by \cite[Theorem 3.1]{M:Y} $f^{\ast}E$ is semistable. In other words, one has the claim.\smallskip

\noindent \textbf{(b) implies (c).} By hypothesis, for any morphism $f\colon C\to X$, $f^{\ast}E$ is semistable, hence\newline
$f^{\ast}\End(E)=f^{\ast}\left(E\otimes E^{\vee}\right)=f^{\ast}E\otimes f^{\ast}E^{\vee}$ is semistable \emph{i.e.} $\End(E)$ is curve semistable. Since $c_1(\End(E))=0$ then $\End(E)$ is nflat (Lemma \ref{lem2.2}.\ref{lem2.2.a}), therefore it is semistable (Lemma \ref{lem2.2}.\ref{lem2.2.b}) and this implies the semistability of $E$. Finally \cite[Propositions 1.2.9 and 1.3]{C:P} proves that $c_2(\End(E))=0$.\smallskip

\noindent \textbf{(c) implies (b).} If $\dim X=1$ there is nothing to prove of course. Let $\dim X=n\geq2$, repeating the reasoning of \cite[Lemma 2.8]{A:C:1} one may assumes $\K\subseteq\C$. Consider the following Cartesian diagram
\begin{displaymath}
\xymatrix{
X_{\C}\ar[r]^{f}\ar[d] & X\ar[d]\\
\spec(\C)\ar[r] & \spec(\K)
},
\end{displaymath}
by hypothesis $\End(E)$ is a degree $0$ semistable Higgs bundle with $c_2(\End(E))=0$, and by the same lemma $f^{\ast}\End(E)$ is a semistable vector bundle over $X_{\C}$ such that $c_2\left(f^{\ast}\End(E)\right)=0$ (cfr. \cite[Theorem 2.7]{A:C:1}). If there exists a morphism $g\colon C\to X$ such that $g^{\ast}E$ is unstable, then by \cite[Lemma 2.8]{A:C:1} $\overline{g}^{\ast}E$ is unstable, where $\overline{g}$ is the base change morphism of $g$ over $\overline{C}=C\times_{\K}\C$. This gives rise a contradiction: $\overline{g}^{\ast}\End(E)$ is the pullback of $f^{\ast}\End(E)$ over $\overline{C}$ and it is unstable, but by \cite[Theorem 1.4]{B:HR} $\overline{g}^{\ast}\End(E)$ is semistable. To avoid this absurd, $E$ is curve semistable.\smallskip

\noindent \textbf{(c) implies (d).} This is trivial.\smallskip

\noindent \textbf{(d) implies (c).} $\End(E)$ is semistable with $\displaystyle\int_Xc_2(\End(E))=0$, while $c_1(\End(E))=0$ of course. Let
\begin{displaymath}
\{0\}=F_0\subsetneqq F_1\subsetneqq\dotsc\subsetneqq F_{m-1}\subsetneqq F_m=\End(E),
\end{displaymath}
be a Jordan-H\"older filtration of $\End(E)$, by \cite[Corollary 6]{L:A} this can be chosen in a way that the quotients $Q_i=F_i/F_{i-1}$ have vanishing Chern classes for any $i\in\{1,\dotsc,m\}$. In other words, $\End(E)$ is an iterating extension of vector bundles with vanishing Chern classes, thus the same statement holds for the Chern classes of $\End(E)$; in particular $c_2(\End(E))=0$.
\end{prf}
\noindent Moreover, the previous theorem is equivalent to the following one.
\begin{theorem}[{cfr. \cite[Corollary 3.2]{B:B:G}}]\label{th4.2}
On a smooth projective variety $X$, the following statements are equivalent.
\begin{enumerate}[a)]
\item\label{th4.2.a} Let $E$ be a curve semistable vector bundle over $X$. Then $E$ is semistable with respect to some polarization $H$ and $c_2\left(\End(E)\right)=0\in\rA^2(X)$.
\item\label{th4.2.b} The Chern classes of any nflat vector bundle over $X$ vanish.
\end{enumerate}
\end{theorem}
\begin{prf}
If \ref{th4.2.a} holds, let $E$ be a nflat vector bundle over $X$. By Lemma \ref{lem2.2}.\ref{lem2.2.b}, $E$ is semistable. Furthermore, applying also \cite[Proposition 2.6.(ii)]{B:GO:1}, $E$ is curve semistable hence ${c_2(E)=c_2\left(\End(E)\right)=0}$. By \cite[Corollary 6]{L:A}, $E$ is extension of vector bundles whose Chern classes vanish, hence the Chern classes of $E$ vanish.\smallskip

\noindent If \ref{th4.2.b} holds, let $E$ be a curve semistable vector bundle over $X$. For any $f\colon C\to X$ one has $f^{\ast}\End(E)\cong\End\left(f^{\ast}E\right)$. As proved in Theorem \ref{th4.1}, $\End(E)$ is curve semistable. Since $\End(E)$ satisfies the hypotheses of Lemma \ref{lem2.2}.\ref{lem2.2.a}, it is nflat hence $c_2\left(\End(E)\right)=0$. By Lemmata \ref{lem2.1}.\ref{lem2.1.b} and \ref{lem2.2}.\ref{lem2.2.b}, $E$ is semistable.
\end{prf}
\begin{remark}
Since Theorem \ref{th4.1} proves that Theorem \ref{th4.2}.\ref{th4.2.a} holds, one has another proof of Theorem \ref{th4.2}.\ref{th4.2.b}. This has been proved by \cite[Proposition 1.3]{C:P} originally.
\end{remark}
\noindent In the Higgs bundles setting, Theorem \ref{th4.1} changes as it follows.
\begin{theorem}\label{th4.3}
Let $\fE=(E,\varphi)$ a rank $r$ Higgs bundle over a smooth projective variety $X$. Consider the following statements:
\begin{enumerate}[a)]
\item $\theta_1(\fE),\dotsc,\theta_{r-1}(\fE)$ are nef;
\item\label{th4.3.b} $\fE$ is curve semistable;
\item\label{th4.3.c} $\fE$ is semistable with respect to some polarization $H$ and $c_2\left(\End(E)\right)=0\in\rA^2(X)$;
\item $\fE$ is semistable with respect to some polarization $H$ and $\displaystyle\int_Xc_2\left(\End(E)\right)\cdot H^{n-2}=0$.
\end{enumerate}
The following implications hold
\begin{displaymath}
\xymatrix{
\text{(a)}\ar@{<=>}[r] & \text{(b)}\ar@{<=}[r] & \text{(c)}\ar@{<=>}[r] & \text{(d)}.
}
\end{displaymath}
\end{theorem}
\noindent It is enough to repeat the proof of Theorem \ref{th4.1}. However if $\fE$ is curve semistable then\newline
$\End(\fE)=(\End(E),\End(\varphi))$ is H-nflat (by Lemmata \ref{lem2.2}.\ref{lem2.2.d} and \ref{lem2.2}.\ref{lem2.2.a}), hence $\fE$ is semistable, but it is unknown whether $c_2(\End(E))=0$.
\begin{remark}
If $\fE$ is semistable and $c_2(\End(E))=0$ then $\fE$ is semistable with respect to any polarization of $X$. This follows from the fact that the semistability of H-nflat Higgs bundles does not depend on the polarization of $X$ (Lemma \ref{lem2.2}.\ref{lem2.2.b}).
\end{remark}
\noindent From all this, it is interesting to study whether the condition $\ref{th4.3}.\ref{th4.3.b}$ implies the condition $\ref{th4.3}.\ref{th4.3.c}$. To simplify the exposition of the corresponding topics, one introduces the following class
\begin{displaymath}
\Delta(E)=\frac{1}{2r}c_2(\End(E))=c_2(E)-\frac{r-1}{2r}c_1(E)^2\in\rA^2(X)\otimes_{\Z}\Q=\rA^2(X)_{\Q}
\end{displaymath}
which is called \emph{discriminant class of} $E$ (cfr. \cite[Theorem 7]{L:A}). I recall the following conjecture.
\begin{conjecture}[\textbf{Bruzzo and Gra\~na Otero Conjecture}]\label{conj1}
Let $\fE$ be a curve semistable Higgs bundle over $X$. Then $\fE$ is semistable with respect to some polarization $H$ and $\Delta(E)=0$.
\end{conjecture}
\noindent To be more precise, the Conjecture \ref{conj1} can be simplified using \cite[Lemma 3.7]{S:CT}, as remarked in \cite{L:LG}, and \cite[Theorem 2.7]{A:C:1}.
\begin{conjecture}\label{conj2}
Let $\fE$ be a curve semistable Higgs bundle over a smooth complex projective surface $X$. Then $\fE$ is semistable with respect to some polarization and $\Delta(E)=0$. 
\end{conjecture}
\noindent In other and simple words, in order to prove Conjecture \ref{conj1} it is enough to study curve semistable Higgs bundles over smooth projective surfaces; and without loss of generality, it is enough also to assume $\K=\C$.
\begin{rem}\label{rem4.1}
The best of our knowledge, the previous conjecture has been proved in the following cases:
\begin{enumerate}[a)]
\item $r=2$, by \cite[Theorems 4.5, 4.8 and 4.9]{B:GO:HR};
\item $X$ has nef tangent bundle, by \cite[Corollary 3.15]{B:LG};
\item $\dim X=2$ and $\kappa(X)\in\{-\infty,0\}$ (the Kodaira dimension of $X$). Indeed, the statement for ruled surfaces follows by \cite[Proposition 3.11]{B:LG}. Since rational surfaces are rationally connected, the statement follows by \cite[Theorem 3.6]{B:LG}. The statement for Abelian surfaces follows by \cite[Corollary 3.8]{B:LG}, the case of K3 surfaces follows by \cite[Theorem 6.4]{B:L:LG} and thus for Enriques surfaces and hyperelliptic surfaces follow by \cite[Proposition 3.12]{B:LG}. All this complete this case;
\item $\dim X=2$, $\kappa(X)=1$ and other technical hypotheses, see \cite[Proposition 5.6]{B:P}; I have generalized this result to all elliptic surfaces (\cite[Theorem 4.18]{A:C:1}), recently;
\item $X$ is a simply-connected Calabi-Yau variety, by \cite[Theorem 4.1]{B:C};
\item if $X$ satisfies the Conjecture \ref{conj1} and $Y$ is a fibred projective variety over $X$ with rationally connected fibres, then $Y$ does the same, by \cite[Proposition 3.11]{B:LG};
\item if $X$ satisfies the Conjecture \ref{conj1} then any finite \'etale quotient $Y$ of $X$ does the same, by \cite[Proposition 3.12]{B:LG};
\item $\fE$ has a Jordan-H\"older filtration whose quotient are H-nflat and have rank at most $2$, by \cite[corollaries 4.3.6 and 4.3.7]{AC:PhD};
\item particular Higgs bundles described in \cite{B:GO:2}. \hfill{$\Diamond$}
\end{enumerate}
\end{rem}
\noindent Furthermore, the previous conjectures are equivalent to a third one (\cite[Corollary 3.2]{B:B:G} and Lemma \ref{lem2.2}.\ref{lem2.2.f}).
\begin{conjecture}\label{conj3}
Let $\fE$ be a stable H-nflat Higgs bundle over a smooth complex projective surface $X$. Then its Chern classes vanish.
\end{conjecture}
\noindent Indeed, to compute Chern classes of $\fE$, one replaces it with graded module $\displaystyle\bigoplus_{i=0}^m\fF_i/\fF_{i+1}$ associated to the filtration \eqref{eq2.1}. Thus and by Remark \ref{rem1.1}.\ref{rem1.1.c}, the second Chern class of $\fE$ is the sum of the second Chern class of the $\fF_i/\fF_{i+1}$'s, which are stable and H-nflat Higgs bundles (Lemma \ref{lem2.2}.\ref{lem2.2.f}).
\begin{remark}
When the Higgs field vanishes, \cite[Proposition 1.3]{C:P} proves the vanishing of Chern classes for nflat vector bundles over smooth projective varieties. This last result has been extended to compact K\"ahler manifolds by \cite[Corollary 1.19]{D:P:S}.
\end{remark}


\begin{thebibliography}{20}
\bibitem{B:B:G} I. Biswas, U. Bruzzo, S. Gurjar - \emph{Higgs bundles and fundamental group schemes}, Adv. Geom. \textbf{19} (2019) 381--388.
\bibitem{B:C} U. Bruzzo, A. Capasso - \emph{Filtrations of numerically flat Higgs bundles and curve semistable Higgs bundles on Calabi-Yau manifolds}, Adv. Geom. \textbf{23} (2023) 215--222.
\bibitem{B:C:GO} U. Bruzzo, A. Capasso, B. Gra\~na Otero - \emph{Positivity for Higgs vector bundles: criteria and applications}, Rev. Mat. Compl. (2025) \url{doi.org/10.1007/s13163-025-00551-7}.
\bibitem{B:GO:1} U. Bruzzo, B. Gra\~na Otero - \emph{Numerically flat Higgs vector bundles}, Commun. Contemp. Math. \textbf{9} (2007) 437--446.
\bibitem{B:GO:2} U. Bruzzo, B. Gra\~na Otero - \emph{Metrics on semistable and numerically effective Higgs bundles}, J. reine ang. Math. \textbf{612} (2007) 59--79.
\bibitem{B:GO:3} U. Bruzzo, B. Gra\~na Otero - \emph{Semistable and numerically effective principal (Higgs) bundles}, Advances in Mathematics \textbf{226} (2011) 3655--3676.
\bibitem{B:GO:HR} U. Bruzzo, B. Gra\~na Otero, D. Hern\'andez Ruip\'erez - \emph{On a conjecture about Higgs bundles and some inequalities}, Mediterr. J. Math. \textbf{20} (2023) article ID 296.
\bibitem{B:HR} U. Bruzzo, D. Hern\'andez Ruip\'erez - \emph{Semistability vs. nefness for (Higgs) vector bundles}, Diff. Geom. Appl. \textbf{24} (2006) 403--416.
\bibitem{B:L:LG} U. Bruzzo, V. Lanza, A. Lo Giudice - \emph{Semistable Higgs bundles on Calabi-Yau manifolds}, Asian J. Math. \textbf{23} (2019) 905--918.
\bibitem{B:LG} U. Bruzzo, A. Lo Giudice - \emph{Restricting Higgs bundles to curves}, Asian J. Math. \textbf{20} (2016) 399--408.
\bibitem{B:P} U. Bruzzo, V. Peragine - \emph{Semistable Higgs bundles on elliptic surfaces}, Adv. Geom. \textbf{22} (2022) 151--169.
\bibitem{C:P} F. Campana, T. Peternell - \emph{Projective manifolds whose tangent bundles are numerically effective}, Math. Ann. \textbf{289} (1991) 169--187.
\bibitem{AC:PhD} A. Capasso, \emph{Positivity Conditions for Higgs Bundles and Applications} Ph.D. Thesis, Universit\`a degli Studi Roma Tre, 2024.
\bibitem{A:C:1} A. Capasso - \emph{Lefschetz principle-type theorems for curve semistable Higgs sheaves and applications to elliptic surfaces} \url{arXiv:2412.00439v3} \textcolor{blue}{\texttt{[math.AG]}}.
\bibitem{A:C:2} A. Capasso - \emph{Curve semistable Higgs bundles and smooth projective varieties whose canonical bundle is ample} \url{arXiv:2508.12062v2} \textcolor{blue}{\texttt{[math.AG]}}.
\bibitem{D:P:S} J.-P. Demailly, T. Peternell, M. Schneider - \emph{Compact complex manifolds with numerically effective tangent bundles}, J. Alg. Geom. \textbf{3} (1994) 295--345.
\bibitem{H:RC} R. C. Hartshorne (1977) \emph{Algebraic Geometry}, Springer.
\bibitem{L:A} A. Langer - \emph{Bogomolov's inequality for Higgs sheaves in positive characteristic}, Inv. Math. \textbf{199} (2015) 889--920.
\bibitem{L:LG} V. Lanza, A. Lo Giudice - \emph{Bruzzo's conjecture}, J. Geom. Phys. \textbf{118} (2017) 181--191.
\bibitem{L:RK} R. K. Lazarsfeld (2004) \emph{Positivity in algebraic geometry. I and II}, Springer.
\bibitem{M:Y} Y. Miyaoka - \emph{The Chern classes and Kodaira dimension of a minimal variety}, Algebraic geometry, Sendai (1985) 449--476. Adv. Stud. Pure Math. \textbf{10}, North-Holland Publishing Co. (1987) Amsterdam.
\bibitem{N:Nak} N. Nakayama - \emph{Normalized tautological divisors of semi-stable vector bundles}, S\=urikaisekikenky\=usho K\=oky\=uroku (1999) 167--173. Kyoto University, Research Institute for Mathematical Sciences.
\bibitem{O:S:S} C. Okonek, M. Schneider, H. Spindler (1988) \emph{Vector Bundles on Complex Projective Spaces. With an Appendix of S. I. Gel'fand}, Birkh\"auser.
\bibitem{S:CT} C. T. Simpson - \emph{Higgs bundles and local systems}, Inst. Hautes \'Etudes Sci. Publ. Math. \textbf{75} (1992) 5--95.
\end{thebibliography}
\end{document}